\documentclass[11pt]{amsart}
\usepackage{graphicx}
\usepackage{color}
\usepackage{setspace}
\usepackage{eufrak}
\usepackage{float}
\usepackage{etex}
\usepackage[all]{xy}
\usepackage{amssymb, amsmath}
\usepackage{ascmac}

\def\qed{\hfill $\Box$}
\def\proof{\noindent {\sl Proof} :\;  }
\def\t{\noindent}

\newcommand{\V}{\mathcal{V}}

\newcommand{\A}{\mathcal{A}}

\newcommand{\Proj}{\mathbb{P}}
\newcommand{\C}{\mathbb{C}}

\newcommand{\R}{\mathbb{R}}

\def\qed{\hfill $\Box$}
\def\proof{\noindent {\sl Proof} :\;  }
\def\t{\noindent}
\def\rd{\partial}

\def\bu{\mbox{\boldmath $u$}}

\def\bp{\mbox{\boldmath $p$}}

\def\b0{\mbox{\boldmath $0$}}

\def\I{{\rm I}}
\def\II{{\rm I\hspace{-.1em}I}}
\def\III{{\rm I\hspace{-.1em}I\hspace{-.1em}I}}
\def\IV{{\rm I\hspace{-.1em}V}}
\def\V{{\rm V}}
\def\VI{{\rm VI}}

\newtheorem{thm}{\bf Theorem}[section]

\newtheorem{prop}[thm]{\bf Proposition}

\newtheorem{rem}[thm]{\bf Remark}

\begin{document}
\title[Family of smooth surfaces in $3$-space]
{Projective classification of jets of surfaces in $3$-space}
%
\author[H.~Sano]{Hiroaki Sano}
\address[H. ~Sano]{Department of Mathematics,
Graduate School of Science,  Hokkaido University,
Sapporo 060-0810, Japan}
\email{s133025@mail.sci.hokudai.ac.jp}
\author[Y.~Kabata]{Yutaro Kabata}
\address[Y. ~Kabata]{Department of Mathematics,
Graduate School of Science,  Hokkaido University,
Sapporo 060-0810, Japan}
\email{kabata@mail.sci.hokudai.ac.jp}
\author[J. L. ~Deolindo Silva]{Jorge Luiz Deolindo Silva}
\address[J. L.~Deolindo-Silva]{Instituto de Ci\^encias Matem\'aticas e de Computa\c{c}\~ao,
Universidade de S\~ao Paulo-S\~ao Carlos -SP, Brazil}
\email{deolindo@icmc.usp.br}

\author[T.~Ohmoto]{Toru Ohmoto}
\address[T.~Ohmoto]{Department of Mathematics,
Faculty of Science,  Hokkaido University,
Sapporo 060-0810, Japan}
\email{ohmoto@math.sci.hokudai.ac.jp}
%
%
%

%
\begin{abstract}
We present
a local classification of smooth surfaces in $\Proj^3$ via
projective transformations in accordance with
singularity types of central projections
with codimension $\le 4$.
We also discuss relations between our classification of Monge forms
and bifurcations of parabolic curves and flecnodal curves. 
\end{abstract}

\maketitle

\section{Introduction}
\subsection{Normal forms}
In projective differential geometry,
local properties of a surface in $3$-space
was intensively investigated by Darboux and others
from the end of the 19th century to the early 20th century, 
while some new interests and ideas
have recently been brought from singularity theory,
generic differential geometry and  applied mathematics
such as computer-vison
\cite{Arnoldbooklet, Arnoldency, BT1, Kabata, Landis, Platonova, Platonova2, Rieger3, UV1, UV2};
for example,
geometric features of a multi-parameter family of surfaces in $3$-space
are of particular interest in application \cite{UV2}. 
In the present paper,
we deal with a classification of jets of surfaces in $\Proj^3$
via projective transformations, which gives a new insight of the classical subject
from a singularity theoretic approach.
We mainly work over $\R$.
Throughout this paper,
we identify $\R^3$ with an open chart of the projective space,
$\{[x:y:z:1]\} \subset \Proj^3$,
and consider germs of surfaces in $\R^3$ at the origin
given by Monge forms $z=f(x, y)$ with $f(0)=0$ and $df(0)=0$.
We say that two germs or jets of surfaces
are {\it projectively equivalent}
if there is a projective transformation on $\Proj^3$ sending one to the other.
A classically well-known fact is  that
at a general hyperbolic point of a surface,
the jet of the Monge form is projectively equivalent to
$$xy+x^3+y^3+\alpha x^4 + \beta y^4 + \cdots$$
where moduli parameters $\alpha, \beta$ are primary projective differential invariants
(those should be compared with the Gaussian/mean curvatures in the Euclidean case and
the Pick invariant in the equi-affine case, cf. \cite{Olver, Tresse, Wil}).
Our aim is to present this kind of expression for degenerate Monge forms.
Most interesting is the parabolic case.

Projective transformations preserving the origin and the $xy$-plane
form a $10$-dimensional subgroup of $PGL(4)$,
and it acts on the space of jets of Monge forms.
We then give an invariant stratification of the jet space of Monge forms 
in accordance with singularity types of central projections, 
and describe the normal form of each stratum up to codimension $4$.
This is a natural extension of results of Platonova
\cite{Platonova} and Landis \cite{Landis}  (see also \cite{Arnoldbooklet, Arnoldency, Goryunov})
 on the classification of jets of Monge forms for a {\it generic} surface,
 which corresponds to the case of codimension $\le 2$ in our list.
 Our main theorem is stated as follows:
\begin{thm}\label{thm1}
Let $M$ be a closed surface
and $U$ an open neighborhood of the origin in $\R^2$ (parameter space).
There is a residual subset $\mathcal{O}$
of the space of $2$-parameter families of smooth embeddings
$\phi: M \times U \to \Proj^3$ equipped with $C^\infty$-topology
so that
each $\phi \in \mathcal{O}$ satisfies that
for arbitrary point $(x_0, u) \in M \times U$
the $p$-jet of Monge form of the surface $M_u (=\phi(M\times u))$ at $x_0$
is projectively equivalent to
one of the normal forms at the origin given in
Tables \ref{main_table1}, \ref{main_table2}, \ref{main_table3}
with codimension $s\; (\le 4)$.
\end{thm}
In Tables,
$\alpha, \beta, \cdots$ are moduli parameters and
$\phi_k$ denotes arbitrary homogeneous polynomials of degree $k$. 
In our stratification, 
each stratum is determined by its {\it projection type}
(i.e. singularity type of projection along the asymptotic line),
that is explained in the next \S 1.2 and \S 2.2. 
This characterization is quite relevant -- e.g. classes $\Pi^p_{\I,k}$ can not be distinguished
by using types of height functions,  the parabolic curves and asymptotic curves,
but can be so by
the difference of singularity types of projection from a special isolated viewpoint. 
The same approach to the classification of surfaces in $\Proj^4$  is considered in \cite{KabataJorge}. 

\begin{table}
{
$$
\begin{array}{l | l | c c | l  }
\mbox{class}  & \mbox{normal form} & \mbox{p} & \mbox{s}
& \mbox{proj. }  \\
\hline
\hline
\Pi^p_{\I,1} & y^2+x^3 + xy^3+\alpha x^4 & 4 & 1 &
\I_2 \, (\I_3)
\\
\Pi^p_{\I,2} & y^2+x^3\pm xy^4+\alpha x^4+\beta y^5 + x^2\phi_3 & 5 & 2 &  \I_2 \, (\I_4)
\\
\Pi^p_{c,1} & y^2+x^2 y + \alpha x^4 \;\;\; (\alpha\not=0, \frac{1}{4}) & 4 & 2 &  \III_2\, (\III_3)
\\
\hline
\Pi^p_{c,2}& y^2+x^2 y + \frac{1}{4} x^4 + \alpha x^5+ y \phi_4 \; (\alpha\not=0) & 5 & 3 &  \III_2 
\\
\Pi^p_{c,4} & y^2+ x^2 y + x^5 + y\phi_4 & 5 & 3 & \IV_1
\\
\Pi^p_{\I,3} & y^2+x^3+ xy^5+\alpha x^4+\phi  & 6 & 3 &  \I_2\,  (\I_{5})\\
\Pi^p_{v,1}& y^2\pm x^4 +\alpha x^3y+\beta x^2y^2 \;\; (\beta\not=\pm\frac{3}{8}\alpha^2) & 4 & 3 & \V_1\,  (\VI)
\\
\Pi^p_{c,3} & y^2+x^2 y + \frac{1}{4} x^4 + y \phi_4  & 5 & 4 &  \III_3\, (\III_{4})
\\
\Pi^p_{c,5} & y^2+ x^2 y \pm x^6+y(\phi_4+\phi_5)    & 6 & 4 &  \IV_2
\\
\Pi^p_{\I,4} &y^2+x^3+\alpha x^4+\phi & 6 & 4 &  \I_2\,  (\I_{6})\\
\Pi^p_{v,2}& y^2 \pm x^4+ \alpha x^3y \pm \frac{3}{8}\alpha^2 x^2y^2  & 4 & 4 & \V_1\,  (\VI_1)\\
\Pi^p_{v,3} & y^2+ x^5 + y(\phi_3+\phi_4)  & 5 & 4 & \V_2\, (\VI_2)
\\
\end{array}
$$
}
\caption{\small  Parabolic case.
In normal forms,
$\phi_k$ denotes homogeneous polynomials of degree $k$
(that is also in Tables \ref{main_table2},  \ref{main_table3}).
For $\Pi^p_{\I,k}\; (k=3,4)$, we put
$\phi=\beta y^5 + \gamma y^6 + x^2(\phi_3 + \phi_4)$ for short.
Double-sign $\pm$ corresponds in same order for each of $\Pi^p_{v,1}$ or $\Pi^p_{v,2}$.
For $\Pi^p_{c,1}$,
if $\alpha=1$, the normal form of $4$-jet should be
$y^2+x^2y+x^4+\beta x^3y$. }
\label{main_table1}
\end{table}

\begin{table}
{
$$
\begin{array}{l | l | c c | l  }
\mbox{class}  & \mbox{normal form} & \mbox{p} & \mbox{s}
& \mbox{proj. }  \\
\hline
\hline
\Pi^h_{3,3} & xy + x^3+y^3 + \alpha x^4+\beta y^4 & 4 & 0 & \II_3 / \II_3
\\
\Pi^h_{3,4} & xy+x^3 + y^4 + \alpha xy^3 & 4 & 1&   \II_3 / \II_4
\\
\Pi^h_{3,5} & xy + x^3+  y^5 +\alpha xy^3  + x\phi_4  & 5 & 2 & \II_3 / \II_5
\\
\Pi^h_{4,4} & xy + x^4\pm y^4+\alpha xy^3+\beta x^3y & 4 & 2 &  \II_4 / \II_4
\\
\hline
\Pi^h_{3,6} & xy+ x^3 + y^6 + \alpha xy^3+x(\phi_4+\phi_5)  & 6 & 3 & \II_3 /\II_6
\\
\Pi^h_{4,5} & xy+x^4+y^5 +\alpha xy^3+\beta x^3y+x\phi_4  & 5 & 3 & \II_4/\II_5
\\
\Pi^h_{3,7} & xy+ x^3 + y^7+\alpha xy^3+x(\phi_4+\phi_5+\phi_6)  & 7 & 4 & \II_3 /\II_7
\\
\Pi^h_{4,6} & xy+ x^4 \pm y^6+\alpha xy^3+\beta x^3y+x(\phi_4+\phi_5)  & 6 & 4 & \II_4/\II_6
\\
\Pi^h_{5,5} & xy + x^5+ y^5+\alpha xy^3+\beta x^3y+ xy\phi_3   & 5 & 4 & \II_5 / \II_5
\end{array}
$$
}
\caption{\small
Hyperbolic case.  
There are two projection types with 
respect to distinct asymptotic lines.  
Here $\II_3$ is of ordinary cusp. 
}
\label{main_table2}
\end{table}

\begin{table}
{
$$
\begin{array}{l | l | c c | l }
\mbox{class}  & \mbox{normal form} & \mbox{p} & \mbox{s}
& \mbox{proj.}  \\
\hline
\hline
\Pi^e & x^2+y^2  & 2 & 0 & \mbox{\footnotesize \rm fold}
\\ \hline 
\Pi^f_{1} & xy^2\pm x^3 + \alpha x^3 y + \beta y^4 & 4  & 3 & {\I_2^\pm}^\dagger
\\
\Pi^f_{2} & xy^2 + x^4 \pm y^4+\alpha x^3 y  & 4 & 4 & {\I_2^-}^\dagger
\\
\end{array}
$$
}
\caption{\small  Umbilical case (elliptic and flat umbilic cases).
For $\Pi^e_{k,j}\; (k,j \ge 4)$, we put $\phi=\beta x^3y + \gamma xy^3$ for short.
$\dagger$: The projection type of $\Pi^f_{k}\; (k=1,2)$ are generically of type $\I_2$, 
see  \S 2.3 for the detail. }
\label{main_table3}
\end{table}

\begin{rem}{\rm
In the case over $\C$,
the elliptic case and the sign difference in
Tables are omitted.
Theorem \ref{thm1} can be restated
in the algebro-geometric context
by mean of a Beltini-type theorem for the linear system of projective surfaces
of degree greater than $p$.
}
\end{rem}

\begin{rem}
{\rm
In the normal forms listed in Tables,
continuous moduli parameters  $\alpha, \beta, \cdots$ and
higher coefficients
must be projective differential invariants in the sense of Sophus Lie;
in fact, there is a similarity between our arguments and
those in a modern theory of differential invariants due to Olver \cite{Olver}.
Besides, those leading parameters may have some particular geometric meanings.
For instance,
at a cusp of Gauss (that is a point of type $\Pi^p_{c,1}$ on the surface),
the coefficient $\alpha$ of $x^4$ in the normal form
coincides with the Uribe-Vargas cross-ratio invariant
defined in \cite{UV1}. Further, we will see that
the coefficient of $x^5$ in the same Monge form
corresponds to the position of a special viewpoint lying on
the asymptotic line (\S 2.2).
Also we discuss in the final section \S 4
a partial connection between our Monge forms
and a topological classification of differential equations (BDE) defining nets of asymptotic curves
(\S 1.3);
indeed $\alpha, \beta, \cdots$ are related to initial moduli parameters of the BDE. 
}
\end{rem}

\subsection{Singularities of projections}
In singularity theory, map-germs are classified up to
diffeomorphism-germs of source and target, that is called
the {\it $\A$-equivalence of map-germs};
we denote by $f \sim_\A g$ for $\A$-equivalent map-germs $f$ and $g$.
Given a surface $M \subset \Proj^3$,
the restriction of linear projection
$\pi_p: \Proj^3-p \to \Proj^2$
from  arbitrary viewpoint $p$ not lying on $M$ is called
the {\it central projection};
for each point $x \in M$, it locally presents a map-germ
$\pi_p|_M: \R^2, 0 \to \R^2, 0$
by taking local coordinates centered at $x$ and $\pi_p(x)$
of source and target, respectively.
Obviously,
if two germs of surfaces are projectively equivalent,
then they admit the same $\A$-types of projections from arbitrary viewpoints.
As a generalization of Arnold-Platonova's result \cite{Platonova2},
 Kabata  studied in \cite{Kabata} central projections of {\it generic families} of surfaces
in the context of $\A$-classification theory of plane-to-plane germs
 \cite{Goryunov, Rieger, Rieger2}.
In particular,
he showed that all $\A$-types of  $\A$-codimension $\le 6$
appear in central projections of generic $2$-parameter families of surfaces\footnote{
It is shown in \cite{Kabata} that
$\A$-types no.$10, 15, 18$ (in Rieger's notation) of codimension $6$ do not appear
in central projections for a generic $1$-parameter family of surfaces.
That is an extension of the fact proven in \cite{Platonova2} that
$\A$-types no.$8, 12, 16$ of codimension $5$ do not appear
in central projections for a generic surface.
}.
We then use part of Kabata's result
and classify jets of degenerate Monge forms
up to projective transformations.
In particular,  Platonova's degenerate classes
$\Pi_*^{(k)}\; (k\ge 1)$ in \cite[Table 2]{Platonova}
break into much finer classes in jet spaces of higher order.

\

\t
{\bf ($\A$-classification)}
Typical singularities are of fold and cusp types,
whose normal forms are given by
$(y, x^2)$ and $(y, xy+x^3)$, respectively.
More complicated $\A$-types of plane-to-plane map-germs
are classified by Rieger \cite{Rieger, Rieger2}.
We follow his notations; $\A$-types with codimension $\le 6$
are named by $1,2, \cdots, 19$ as in  \cite[Table 1 (p.352)]{Rieger}.
Some $\A$-types are combined into a single topological-$\A$-type,
that was studied in detail  in \cite{Rieger2};
they are characterized by some {\it specified $p$-jets} listed
in Table \ref{Riegerlist_top} below, and called
{\it equisingularity types}  $\I_k, \cdots, \VI$ (see \cite{Kabata}).
For our convenience,
provisionally we add two types of $4$-jets
with codimension greater than $6$, which are adjacent to $\VI$ ($c\not=0$):
$$\VI_1:\; (y, x^4+d xy^3), \quad \VI_2: (y, c x^2y^2+d xy^3).$$

\begin{table}
$$
\begin{array}{c | c |l | c | c }
\mbox{type} & \mbox{$\A$-type}   & \mbox{normal forms of $p$-jets} & p& \mbox{codim.}\\
\hline
\I_k  & 4_k^\pm & (y, x^3\pm xy^k) &k+1& k+1 \\
\II_k  & 5-10 & (y, xy+x^k) &k& k-1\\
\III_k& 11_{2k+1} & (y, x^2y+x^4+x^{2k+1}) &2k+1& k+2 \\
\IV_1 & 12, 13, 14  & (y, x^2y+x^5) &5&5\\
\IV_2 & 15  & (y, x^2y+x^6)&6&6 \\
\V_1 & 16, 17  & (y, xy^2+x^4) &4&5 \\
\V_2 & 18  &  (y, xy^2+ d x^3y) &4&6 \\
\VI & 19 & (y, x^4+ c x^2y^2+d xy^3) &4&6
\end{array}
$$
\caption{\small Equisingularity types of plane-to-plane germs
with $\A$-codimension $\le 6$  \cite{Rieger2}}
\label{Riegerlist_top}
\end{table}

\

\t
{\bf (Projection type)}
The projection type in Tables
means the $\A$-type of the central projection from
{\it arbitrary} viewpoint lying on the asymptotic line.
Here {\it the entire Monge form is assumed to have
appropriately generic higher terms of order $>p$
added to the prescribed normal forms of $p$-jet in Tables. }

For instance, look at the class $\Pi^p_{c,1}$;
the stratum has codimension $2$ in the jet space, thus
the class appears at some isolated point on a generic surface,
that is classically called an (ordinary) {\it cusp of Gauss} point.
In Table \ref{main_table1}
the corresponding  projection type is written as  ``$\III_2\, (\III_3)$".
That means that the central projection from
almost all viewpoints lying on the asymptotic line
has the Gulls singularity $\III_2$,
while the projection from some exceptional isolated viewpoint
is of type  $\III_3$ worse than Gulls type.
If $\alpha=0$, we have
the degenerate class $\Pi^p_{c,4}$; it has codimension $3$, so
it appears generically in $1$-parameter family of surfaces
with the projection type $\IV_1$.
In case that $\alpha=\frac{1}{4}$,
there are two degenerate classes  $\Pi^p_{c,2}$ and $\Pi^p_{c,3}$
according to whether the term $x^5$ in the $5$-jet remains or not.
The latter type appears generically in a $2$-parameter family of surfaces,
and the projection from any viewpoints
on the asymptotic line is of type $\III_3$ or worse. See \S 2.2 for the detail.

\subsection{Net of asymptotic curves}
At each hyperbolic point of a surface in $\Proj^3$,
there are exactly two asymptotic lines (lines tangent to the surface
with more than $2$-point contact); they are invariants
under projective transformations.
The integral curves, called {\em asymptotic curves}, form
a pair of foliations on the hyperbolic domain,
which is classically named by the {\em net of asymptotic curves}.
The {\it parabolic curve} is
the locus of singular points of asymptotic curves,
and in fact it is the locus of points on the surface
whose Monge form is of type $\Pi^p_{\I,1}$ or more degenerate ones. 
The {\it flecnodal curve} 
is the locus of inflection points of asymptotic curves; 
it is actually the closure of the locus of class $\Pi^h_{3,4}$. 
The parabolic and flecnodal curves meet each other tangentially 
at a cusp of Gauss, i.e. a point of class $\Pi^p_{c,k}\; (k\ge 1)$. 

The net of asymptotic curves is defined
by a {\it binary differential equation} (BDE).
In a general setting,
Davydov \cite{Davydov, Davydov2} and Bruce-Tari \cite{BT1, BT2, BFT, Tari1, Tari2} has presented
the topological classification of (families of) BDE. 
We then compare our classification of parabolic Monge forms with 
degenerate parabolic points arising in the general classification of BDE
(Propositions \ref{BDE0}, \ref{BDE1}, \ref{BDE2}).  
In particular, we will see that
the flat umbilic class $\Pi^f_{2}$ corresponds to a type of 
generic $3$-parameter family of BDE  studied in Oliver \cite{Oliver}. 
Furthermore, 
we also compare our classification of hyperbolic Monge forms with bifurcations of flecnodal curves. 
In part of his dissertation \cite{UV2} 
Uribe-Vargas presented a complete list 
of generic $1$-parameter bifurcations of flecnodal curves 
via a geometric approach using the dual surfaces. 
We give a precise characterization 
of moduli parameters in our corresponding Monge form 
for each type of Uribe-Vargas' classification.

\subsection*{Acknowledgement}
The authors would like to thank Takashi Nishimira and Farid Tari
for organizing the JSPS-CAPES international cooperation project in 2014-2015.
In fact, the second and third authors are supported by the project
for their stays in ICMC-USP and Hokkaido University, respectively.  
The authors also thank Ricardo Uribe-Vargas for letting us know of his paper \cite{UV2} and 
for his valuable comments. 
They are partly supported by JSPS grants no.24340007 and 15K13452.

\section{Jets of Monge forms}

\subsection{Central projection}
Assume that $\R^4$ is equipped with the standard inner product.
Let $p=[\bp]\in \Proj^3$, called a {\it viewpoint}, and
let $W_p \subset \R^4$ denote
the orthogonal complement to the vector $\bp\in \R^4$.
The {\it central projection} $\pi_p$ is the map from $\Proj^3-\{p\}$
to the projectivization $\Proj(W_p)$
given by
$$\pi_p([\bu])=\left[\bu-\frac{(\bu \cdot \bp)\bp}{||\bp||^2}\right] \in \Proj(W_p) \subset \Proj^3.$$
Restrict $\pi_p$ to the open set $\R^3 \subset \Proj^3$.
For $\bp=(a,b,c,1) \in \R^4$,
set
$$
A=\left(\begin{array}{rrrr}
1 & 0 & 0 & -a \\
0 & 1 & 0 & -b \\
0 & 0 & 1 & -c \\
-a & -b & -c & -1
\end{array}\right)
$$
and $\Phi_A: \Proj^3 \to \Proj^3$ the projective transformation
defined by $A$.
Obviously, $\Phi_A(p)=0 \in \R^3$ and $\Phi_A(W_p)=\Proj^2$.
We identify $\pi_p$ with
$$\Phi_A\circ \pi_p: \R^3-\{p\} \to \Proj^2, \quad (x,y,z)\mapsto [x-a:y-b:z-c].$$
If the viewpoint is at infinity, i.e. $\bp=(p,0)$ with $p=(a,b,c) \in \R^3$,
then the projection is given by for $\bu=(u, 1) \in \R^4$
$$\pi_p([\bu])=\left[\bu-\frac{(\bu \cdot \bp)\bp}{||\bp||^2}\right]=\left[u-\frac{(u \cdot p) p}{||p||^2}:1\right] \;\; \in \Proj(W_p).$$
Hence it induces the {\it orthogonal  projection}  (or parallel projection)
in $\R^3$ along the line generated by the vector $p$;
if $a\not=0$ and $v=b/a$, $w=c/a$, then we have
by a linear transform on target $\R^2$
$$\pi_p: \R^3 \to \R^2, \quad \pi_p(x,y,z)=(y-vx, z-wx).$$

Let $M$ be a surface in $\R^3 \; (\subset \Proj^3)$
around the origin with the Monge form
$$z=f(x,y)=\sum_{i+j\ge 2} c_{ij}x^iy^j.$$
Take a viewpoint $p=(a,b,c)$ with $a\not=0$, and
then the central projection from $p$ is locally written by
$$\varphi_{p,f}=\pi_p|_M: M \to \R^2, \;\;\;
(x,y)\mapsto \left(\frac{y-b}{x-a},\frac{f(x,y)-c}{x-a}\right)=(X,Y)
$$
using local coordinates $(x,y)$ of $M$ and $[1:X:Y]$ of $\Proj^2$.
If $p$ is chosen to be at infinity in $\Proj^3$,
then
$$\varphi_{p,f}(x,y)=\left(y - v x,f(x,y) - w x\right).$$

Given a surface $M$ and a point $x \in M$,
we are interested in the singularity types
(as plane-to-plane map-germs) of central projections
$\varphi_{p,f}:M, x\to \R^2,0$
from arbitrary viewpoint $p \in \Proj^3-M$.
We say two germs $\eta_1, \eta_2: \R^2, 0 \to \R^2, 0$
are {\it $\A$-equivalent} if
there are diffeomorphism-germs $\sigma, \tau$ of source and target $\R^2$ at the origins
such that
$\eta_2=\tau \circ \eta_1 \circ \sigma^{-1}$.

\subsection{Criteria} \label{criteria}
Now, assume that the $x$-axis is an asymptotic line of $M$ at the origin
and $p$ lies on it, i.e. $b=c=0$.
Then, taking new coordinate $\bar{x}=x$ and $\bar{y}=y(x-a)^{-1}$, the projection is of the form
$$ \varphi_{p,f}=\left(\bar{y}, \frac{f(\bar{x}, \bar{y}(\bar{x}-a))}{\bar{x}-a}\right).$$
According to Rieger's classification \cite{Rieger}
and useful criteria in \cite{Kabata} for $\A$-types of plane-to-plane germs,
one can determine
local types of singularities of $\varphi_{p,f}$.
In fact, for each $\A$-type in Rieger's list,
Kabata explicitly described in \cite{Kabata}
the necessary and sufficient condition on coefficients $c_{ij}$
so that the germ of $\varphi_{p,f}$ at the origin is $\A$-equivalent to the type.
Elliptic case and umbilical case are easy, so we omit them.
For the hyperbolic and parabolic cases,
we may assume that $j^2f=xy$ and $j^2f=y^2$, respectively,
by some linear transformation.
In Tables \ref{criteria_h} and \ref{criteria_p},
the middle column stands for Kabata's closed condition on $c_{ij}$
and the right is the corresponding projection type,
i.e. the $\A$-type of $\varphi_{p, f}$.
The condition defines an invariant stratum named by class $\Pi^*_{*,*}$ (left).
Notice that some open condition on $c_{ij}$ is implicitly imposed
for each class in the list, e.g., in Table \ref{criteria_h},
we understand that $c_{40}\not=0$ for $\Pi^h_{3,4}$
to distinguish the type from others.

\begin{table}
$$
\begin{array}{l  l  l}
&\mbox{\footnotesize cond.} &
\mbox{\footnotesize type }\\
\hline
\Pi^h_{3,4} &
 c_{30}=0 & \II_4\\
\Pi^h_{3,5} &
c_{30}=c_{40}=0 & \II_5\\
 \Pi^h_{3,6} &
c_{30}=c_{40}=c_{50}=0 & \II_6   \\
\Pi^h_{3,7} &
c_{30}=c_{40}=c_{50}=c_{60}=0 & \II_7\\
\end{array}
$$
\caption{\small
Criteria of ($\A$-) equiangularity types
for projections at a hyperbolic point \cite{Kabata}.
It is assumed that $c_{11}\not=0$ and $c_{20}=c_{02}=0$.
If $c_{30}, c_{03}\not=0$, then it is a general hyperbolic point ($\Pi^h_{3,3}$),
where the projection is of type cusp.
The flecnodal curve consists of points of type
$\Pi^h_{3,k}\; (k\ge 4)$,
and its self-intersection points are of type $\Pi^h_{4,4}$ etc.
 }
\label{criteria_h}
\end{table}

\begin{table}
$$
\begin{array}{l  l  l}
&\mbox{\footnotesize cond.} &
\mbox{\footnotesize type}\\
\hline
\Pi^p_{\I,1} &
c_{20}=0 & \I_2\, (\I_3)\\
\Pi^p_{\I,2} &
c_{20}=P=0 &  \I_2\, (\I_4)\\
\Pi^p_{\I,3} &
c_{20}=P=Q=0 &  \I_2\, (\I_{5})\\
\Pi^p_{\I,4} &
c_{20}=P=Q=R=0 &  \I_2\, (\I_6)\\
\hline
\Pi^p_{c,1} &
c_{20}=c_{30}=0 & \III_2\, (\III_{3})\\
\Pi^p_{c,2}&
c_{20}=c_{30}=B=0 & \III_2 \\
\Pi^p_{c,3}&
c_{20}=c_{30}=A=B=0 & \III_3 \, (\III_{4})\\
\Pi^p_{c,4} &
c_{20}=c_{30}=c_{40}=0 & \IV_1\\
\Pi^p_{c,5} &
c_{20}=c_{30}=c_{40}=c_{50}=0 & \IV_2\\
\hline
\Pi^p_{v,1} &
c_{20}=c_{30}=c_{21}=0 & \V_1\, (\VI)\\
\Pi^p_{v,2} &
c_{20}=c_{30}=c_{21}=S=0 & \V_1\, (\VI_1)\\
\Pi^p_{v,3} &
c_{20}=c_{30}=c_{21}=c_{40}=0 & \V_2\,  (\VI_2)\\
\end{array}
$$
\caption{\small
Criteria of ($\A$-) equiangularity types
for projections at a parabolic point
\cite{Kabata}.
It is assumed that $c_{11}=0$ and $c_{02}\not=0$, and
$A, B, P, \cdots, S$ are some polynomials in $c_{ij}$.
A point of type $\Pi^p_{c,k}\; (k\ge 1)$ is an ordinary/degnerate cusp of Gauss point.}
 \label{criteria_p}
\end{table}

In Table \ref{criteria_p} of the parabolic case,
$A, B, P, \cdots$ are polynomials in $c_{ij}$ given as follows.
They are defined by Kabata's criteria for
the appearance of singularity types
$\I_k$, $\III_k$ and $\VI$ in central projections from
{\it special isolated viewpoints}.
We explain it below
(see also \cite{Kabata});
it is always assumed that $c_{11}=0$, $c_{02}\not=0$
and $p=(a,0,0)$ with $a\not=0$ (possibly $a=\infty$).

\

\t
({\bf $\I_k$})
\;\;
$j^3\varphi_{p,f}\sim_{\A^3} (y,x^3+ \delta xy^2)$ for some $\delta$
if and only if  $c_{20}=0$ and $c_{30}\not=0$.   \\
Here $\delta$ is written by a non-zero scalar multiple of $C a + D$ with
$$C=c_{21}^2-3c_{30}c_{12}, \quad D=3c_{02}c_{30}\;\; (\not=0)$$
(if  $p$ is at the infinity, that is, $\varphi_{p, f}$ is the parallel projection,
$\delta=C/c_{30}^2$).
The type $\I_2$ (lips/beaks) is just the case of $\delta \not=0$.
Therefore,
the projection is of type $\I_2$
from almost all viewpoints lying on the asymptotic line ($x$-axis),
i.e. $C a + D \not=0$,
while there is a unique viewpoint $p=(-D/C, 0, 0)$
(or the infinity if $C=0$)
so that the projection is of type $\I_{3}$ or worse.
Let $p$ be the special viewpoint.
Then the degenerate projection is written
by $\varphi_{p,f}=(y, x^3+\sum_{i+j\ge 4} a_{ij}x^iy^j)$
in some local coordinates; we put
$$\textstyle P=a_{13}, \;\; Q=a_{14}-\frac{1}{3}a_{22}^2, \;\;
R=a_{15}-\frac{2}{3}a_{23}a_{22}+\frac{1}{3}a_{31}a_{22}^2.$$
It is shown in \cite[Prop. 3.8]{Kabata} that

\begin{center}
\begin{tabular}{clcl}
$\bullet$ & $P\not=0$ & $\Longleftrightarrow$ & $\varphi_{p,f} \sim_\A \I_3$ (Goose type);
\\
$\bullet$ & $P=0, \;\; Q\not=0$ & $\Longleftrightarrow$ & $\varphi_{p,f} \sim_\A \I_4$;
\\
$\bullet$ & $P=Q=0, \;\; R\not=0$ & $\Longleftrightarrow$& $\varphi_{p,f} \sim_\A \I_5$;
\\
$\bullet$ & $P=Q=R=0$ & $\Longleftrightarrow$ & $\varphi_{p,f} \sim_\A \I_{\ge 6}$.
\end{tabular}
\end{center}
Notice that $a_{ij}$ are polynomials in $c_{ij}$, hence $P, Q, R$ are so.
These conditions define our classes $\Pi^p_{\I, k}\; (1\le k \le 4)$.

Note that $C=D=0$ if and only if $c_{30}=c_{21}=0$,
that is the case of $\Pi^p_{v,k}$,
where the projection type becomes to be more degenerate, see $(\VI)$ below.

\

\t
({\bf $\III_k$}) \;\;
$j^5\varphi_{p,f} \sim_{\A^5} (y, x^2y+x^4+\delta x^5)$
if and only if
$$
c_{20}=c_{30}=0, \;\;
 c_{40}\not=0, \;\; c_{21}\not=0.$$
Here $\delta=0$ if and only if $Aa+B=0$ with
$$A=c_{21}^2c_{50}+4c_{12}c_{40}^2-2c_{21}c_{31}c_{40}, \quad
B=c_{21}^2c_{40}-4c_{02}c_{40}^2$$
(if $p$ is at the infinity, $\delta$ is a non-zero multiple of $A$).
By the definition,
the germ is of type $\III_2$ (Gulls) if and only if $\delta\not=0$,
and the condition $\delta=0$ determines the special viewpoint
for type $\III_3$ or worse. There are three cases:
\begin{itemize}
\item $\Pi^p_{c,1}$:
In case of $B\not=0$, the projection is of type $\III_2$
from almost all viewpoints on the asymptotic line
except for
a unique viewpoint $p=(-B/A, 0, 0)$ (or at the infinity, if $A=0$);
\item $\Pi^p_{c,2}$:
In case of $B=0$ and $A\not=0$,
there is no special viewpoint --
the projection is of type $\III_2$
from {\it any} viewpoint on the line;
\item $\Pi^p_{c,3}$:
In case of $A=B=0$,
the type $\III_{\ge 3}$ is observed from any viewpoint on the line.
\end{itemize}

\

\t
({\bf $\VI$})  \;\;
 $j^4\varphi_{p,f}\sim_{\A^4}(y, \delta xy^2+ c_{40} x^4+c x^2y^2+d xy^3)$
if and only if
$$c_{20}=c_{30}= c_{21}=0,$$
where $\delta=c_{12}a - c_{02}$ and
 $c$ is a non-zero scalar multiple of
$$S=3 c_{02} c_{31}^2 + 8 c_{40}(c_{12}^2  - c_{02} c_{22}).$$
In this case,
for general viewpoint $p=(a,0,0)$ with $\delta\not=0$,
if $c_{40}\not=0$, the projection $\varphi_{p,f}$ is of type $\V_1:(y, xy^2+x^4)$,
and if $c_{40}=0$, it is of type $\VI$.
Let  $\delta=0$, that is,
$p=(c_{02}/c_{12}, 0, 0)$ (possibly at the infinity) be the special viewpoint.
Then, according to our convention of $\A$-types of codimension $\ge 6$,
\begin{itemize}
\item $\Pi^p_{v,1}$:
If $\delta=0$, $S\not=0$ and $c_{40}\not=0$, $\varphi_{p,f}$ is of type $\VI$;
\item $\Pi^p_{v,2}$:
If $\delta=S=0$ and $c_{40}\not=0$, $\varphi_{p,f}$ is of type $\VI_1$;
\item $\Pi^p_{v,3}$:
If $\delta=c_{40}=0$,  $\varphi_{p,f}$ is of type $\VI_2$.
\end{itemize}
The difference between $\Pi^p_{v,1}$ and $\Pi^p_{v,2}$ is clear geometrically.
A simple computation shows that in the above case (i.e. $c_{11}=c_{20}=c_{30}= c_{21}=0$),
the parabolic curve of the surface near the origin is defined by
$$6 c_{02} c_{40} x^2 + 3 c_{02} c_{31} x y + (c_{02} c_{22}- c_{12}^2) y^2+ h.o.t=0,$$
thus $S$ coincides with the Hessian determinant.
Hence, for $\Pi^p_{v,1}$ (resp. $\Pi^p_{v,2}$),
the parabolic curve has a Morse singularity (resp. a cusp singularity).
For $\Pi^p_{v,3}$, the parabolic curve has a Morse singularity, but the hight function is
more degenerate ($A_4$-singularity) than
those for $\Pi^p_{v,1}$ and $\Pi^p_{v,2}$ ($A_3$-singularity).

\begin{rem}{\rm
If we take attention in the difference of $\A$-types
belonging to the same equisingularity type,
it leads to a more finer stratification.
For instance,
look at a surface of type $\Pi^h_{3,5}$.
We observe Butterfly $6: (y, xy+x^5+x^7)$ from general viewpoints
and Elder Butterfly $7: (y, xy+x^5)$ from at most two special viewpoints,
where these $\A$-types are of the same equiangularity type $\II_5$ --
indeed they are distinguished in Platonova's list \cite{Platonova}.
More precisely,
let $c_{11}\not=0$ and $c_{20}=c_{02}=0$,
then $\Pi^h_{3,5}$ is characterized by $c_{30}=c_{40}=0$,
and there is a quadric equation $A_1a^2+A_2a+A_3=0$,
where $A_i$ are polynomials of $c_{ij}$,
so that the solution $a$ gives
the special position of points for viewing Elder Butterfly.
Obviously,
the number of such viewpoints  ($0, 1$ or $2$) is
a projective invariant, so $\Pi^h_{3,5}$ can be divided
into two open subsets and a closed subset of codimension $1$
in the stratum $\Pi^h_{3,5}$.
}
\end{rem}

\begin{rem}{\rm
It would be interesting to take
the same approach for studying projections of singular surfaces in $3$-space,
e.g. the image of  a certainly `nice' map-germ  $\R^2,0 \to \R^3,0$.
Indeed West \cite{West} studied parallel projections of
a generic crosscap up to affine transformations of target $3$-space
(see  \cite{YKO, YKO2} for non-generic crosscaps).
}
\end{rem}

\subsection{Flat umbilic points} \label{umbilic points}
For a flat umbilic point, any tangent line is asymptotic, thus
for any viewpoint lying on the tangent plane,
the projection type should be considered.
Let $z=f(x,y)$ be the flat umbilic Monge form at the origin and $p=(a,b,0)$,  
and then the following are verified by Kabata's criteria \cite{Kabata} 
(the argument below is also true for the case that $p$ is at the infinity). 

First, consider the class $\Pi^f_{1}$ with codimension $3$.
The sign in the normal form makes a difference.
Let $f=xy^2+x^3+\cdots$.
Then 
\begin{itemize}
\item
for general tangent lines ($a\not=\pm b$), $\varphi_{p,f}$ is of type $\I_2^\pm$;
\item
There are two exceptional tangent lines, $a=\pm b$,
where $\varphi_{p,f}$ is of type $\I_3$ from almost all viewpoints on the line
and of type $\I_4$ from some special points
(it is of type $\I_5$ from special points for the case of codimension $4$).
\end{itemize}
Let $f=xy^2-x^3+\cdots$. Then
\begin{itemize}
\item
for general tangent lines ($a\not=\pm b$), $\varphi_{p,f}$ is of type $\I_2^-$ (beaks only);
\item
If $a=\pm b$, $\varphi_{p,f}$ is of type $\III_2$ from almost all viewpoints on the line
and of type $\III_3$ from some special points
(type $\III_4$ from special points for the case of codimension $4$),
besides, it also happen in codimension $4$ that
$\varphi_{p,f}$ is of type $\IV_1$ from all points on the line.
\end{itemize}
The case of $\Pi^f_{2}$ is easier. Let $f=xy^2+x^4+\cdots$.
From almost all viewpoint on the tangent plane ($b\not=0$), only $\I_2^-$ is observed.
On the line $b=0$, only $\V_1$ is observed.

\section{Proof of Theorem \ref{thm1}}

\subsection{Projective equivalence}
Let $0 \in \R^3 \subset \Proj^3$.
A projective transformation on $\Proj^3$ preserving the origin and the $xy$-plane
defines a diffeomorphism-germ $\Psi: \R^3, 0 \to \R^3,0$ of the form
$$ \Psi(x,y,z)=\left(\frac{q_1(x,y,z)}{p(x,y,z)},
\frac{q_2(x,y,z)}{p(x,y,z)},
\frac{q_3(x,y,z)}{p(x,y,z)}\right)$$
where
$$q_1=u_1 x + u_2 y + u_3 z, \;\; q_2=v_1 x + v_2 y+v_3z, \;\;q_3=cz,$$
$$p=1+w_1 x + w_2 y + w_3 z.$$
Note that there are 10 independent parameters $u_1, \cdots, w_3$.

For Monge forms $z=f(x,y)$ and $z=g(x,y)$,
we denote $j^kf \sim j^kg$
if $k$-jets of these surfaces at the origin are projectively equivalent.
That is equivalent to that there is some projective transformation
$\Psi$ of the above form so that
\begin{equation}\label{eq1}
F(x,y,g(x,y))=o(k)
\end{equation}
where $F(x,y,z)=q_3/p-f(q_1/p, \,q_2/p)$
and $o$ means Landau's symbol.

We shall simplify jets of the Monge forms via projective transformations.
Let $f(x,y)=\sum_{i+j\ge 2} c_{ij}x^iy^j$.
For a simplified form $g(x,y)$,
the equation (\ref{eq1}) yields
a number of algebraic equations in $u_1, \cdots, w_3$,
and then our task is to find at least one solution in terms of $c_{ij}$
(we use the software {\it Mathematica}  for the computation).

\subsection{Parabolic case}
Suppose $c_{11}=c_{20}=0$ and $c_{02}\not=0$.

\

\t
($\Pi^p_{\I,k}$):
Let $c_{30}\not=0$,
then we may assume that $j^3f=y^2+x^3$. Further,
$$j^4f \sim y^2+x^3+c_{40} x^4+c_{13}xy^3$$
by $\Psi$ with
$$\textstyle q_1=x-\frac{1}{3}c_{22}z, \;\;
q_2=y+\frac{1}{2}c_{31}z, \;\; q_3=z,$$
$$\textstyle p=1+c_{31}y+ (c_{04}+\frac{1}{4}c_{31}^2)z.$$
Let $j^4f$ be of the form in the right hand side.
Then $c_{21}=c_{12}=0$ in this form,
thus $C=0$ and $D\not=0$
where $C$ and $D$ are polynomials in previous section,
so the special viewpoint $p$ is at the infinity.
Therefore $\varphi_{p, f}$ is the parallel projection
and
$$\varphi_{p, f} \sim_\A (y, x^3+c_{40} x^4+c_{13}xy^3+h.o.t.).$$
For $c_{22}=0$, we have
$P=c_{13}$, $Q=c_{14}$ and $R=c_{15}$.
Thus the conditions on $P$, $Q$ and $R$
lead to the normal forms of $\Pi^p_{\I,k}\; (1 \le k \le 4)$.

\

\t
($\Pi^p_{c,k}$):
Let $c_{30}=0$ and $c_{21}\not=0$,
then we may assume $j^3f=y^2+x^2y$.
Then
$$j^4f \sim y^2+x^2y+c_{40} x^4+c_{31} x^3y$$
by $\Psi$ with
$$\textstyle q_1=x-\frac{1}{2}c_{13} z,\;\;
q_2=y+c_{22}z, \;\;
q_3=z,$$
$$\textstyle p=1+2c_{22}y+(c_{04}+c_{22}^2)z.$$
Furthermore,  the term $x^3y$ in the right hand side is killed\footnote{
As an exception, if $c_{40}=1$,  we see by the same manner that
$x^3y$ can not be killed via any projective transformations,
so the normal form is $y^2+x^2y+x^4+c_{31}x^3y$. }
by a new $\Psi$ with
$$\textstyle
q_1=x -\frac{c_{31}}{4(1-c_{40})}y+\frac{c_{31^2}}{32(1-c_{40})^2}z,
\; \;
q_2=y+\frac{3c_{31}^2}{8(1-c_{40})}z,
\; \; q_3=z,$$
$$\textstyle
p=1 -\frac{c_{31}}{2(1-c_{40})}y-\frac{c_{31}^2(12c_{40}-13)}{16(1-c_{40})^2}y-\frac{c_{31}^4(36c_{40}-37)}{256(1-c_{40})^3}z, $$
then we have the normal form of $4$-jet:
$$\Pi^p_{c,1}: \;\; y^2+x^2y+c_{40} x^4.$$
As for the $5$-jet,  we define  $\Pi^p_{c,2}$
by adding one more equation $B=0$
and  $\Pi^p_{c,3}$
by $A=B=0$, where $A, B$ are given in the previous subsection.
Notice that the conditions are invariant under projective transformations.
Therefore we may assume that $f=y^2+x^2y+c_{40} x^4+h.o.t.$,
and then $B=0$ implies $c_{40}=\frac{1}{4}$ and $A=c_{50}$.
Hence  we have the normal form for both classes:
$$\textstyle y^2+x^2y+\frac{1}{4}x^4 + c_{50}x^5 + y\phi_4 \;
(=(y+x^2)^2+ c_{50}x^5 + y\phi_4).$$

Let $c_{30}=c_{40}=0$ and $c_{21}\not=0$, then
we have  $\Pi^p_{c,4}$ or $\Pi^p_{c,5}$
according to whether $c_{50}\not=0$ or $0$.

\

\t
($\Pi^p_{v,k}$):
Let $c_{30}=c_{21}=0$, then
we may assume that $j^3f=y^2$.
Further,
if $c_{40}\not=0$,  we may assume $c_{40}=\pm 1$.
Then we have the normal form
$$\Pi^p_{v,1}: \;\; y^2\pm x^4+ \alpha x^3y + \beta x^2y^2$$
for some $\alpha, \beta$ by $\Psi$ with
$q_1=x+u_1y$,
$q_2=y$,
$q_3=z$ and
$$p=1+(\pm 1+c_{13}u_1+c_{22}u_1^2+c_{31}u_1^3 \pm u_1^4)z,$$
where $u_1$ satisfies $c_{13}+2c_{22}u_1+3c_{31}u_1^2 \pm 4u_1^3=0$.
As a degenerate type, $\Pi^p_{v,2}$ is defined by adding $S=0$,
where $S$ is in $(\VI)$ of the previous section.
For the above form of $\Pi^p_{v,1}$,
$S=3 c_{02} c_{31}^2 + 8 c_{40}(c_{12}^2  - c_{02} c_{22})=3\alpha^2\pm 8\beta$,
hence $\Pi^p_{v,2}$ is the case that $\beta=\pm \frac{3}{8}\alpha^2$.

If $c_{40}=0$, we have $j^4f=y^2+y\phi_3$, that defines the class $\Pi^p_{v,3}$.

\subsection{Hyperbolic case}
Notice that  normal forms in Table \ref{main_table2} are obtained 
by interchanging variables $x, y$ used in the following computation. 
Suppose $c_{11}\not=0$ and $c_{20}=c_{02}=0$.

\

\t
($\Pi^h_{3,k}$):
 Let $c_{30}, c_{03}\not=0$, then we may assume
$j^3f=xy+x^3+y^3$.
For the $4$-jet, we obtain the normal form
$$\Pi^h_{3,3}: \; xy+x^3+y^3+\alpha x^4+\beta y^4$$
by $\Psi$ with
$q_1=x-\frac{c_{31}}{2}z$, $q_2=y-\frac{c_{13}}{2}z$, $q_3=z$ and 
$p=1-\frac{c_{13}}{2}x-\frac{c_{31}}{2}y+(c_{22}+\frac{c_{13}c_{31}}{4})z$.

Let $c_{30}=0$ and $c_{03}\not=0$.
By a linear transformation, we may assume $j^3f=xy+y^3$.
Then
$$j^4f \sim xy+y^3+c_{40}x^4+c_{31} x^3y$$
by $\Psi$ with
$$\textstyle q_1=x+c_{04}z, \;\; q_2=y-\frac{1}{2}c_{13}z, \;\; q_3=z, $$
$$\textstyle p=1-\frac{1}{2}c_{13}x+c_{04}y + (c_{22}-\frac{1}{2}c_{13}c_{04})z.$$
If $c_{40}\not=0$,
we have
$$\Pi^h_{3,4}: \; xy+y^3+x^4+ \alpha x^3y.$$
If $c_{30}=c_{40}=0$,
we have the normal form of $\Pi^h_{3,k}\; (k=5,6,7)$ as in Table \ref{main_table1}
according to whether $c_{50}$ and/or $c_{60}$ vanish or not.

\

\t
($\Pi^h_{4,k}$):
Let $c_{30}=c_{03}=0$. Then
$$j^4f \sim xy+c_{40}x^4+\alpha x^3y + \beta xy^3+c_{04} y^4$$
by $\Psi$ with
$$q_1=x-c_{12}z, \;\; q_2=y-c_{21}z, \;\; q_3=z,$$
$$p=1+ (c_{22}-3c_{12}c_{21})z, \;\;
\mbox{with}\;\;
\alpha=c_{31}-c_{21}^2, \;\; \beta=c_{13}-c_{12}^2.$$
This is the case that
each of two asymptotic lines ($x$-axis and $y$-axis)
has degenerate tangencies with the surface;
the normal forms of $\Pi^h_{4,4}$, $\Pi^h_{4,5}$, $\Pi^h_{4,6}$  and $\Pi^h_{5,5}$
are immediately obtained.

\subsection{Flat umbilic case}
Suppose $c_{11}=c_{20}=c_{02}=0$.
Cubic forms of variables $x, y$ are
classified into $GL(2)$-orbits of
$xy^2\pm x^3$, $xy^2$, $y^3$ and $0$.
Let $D$ be the discriminant of the cubic form:
$$D=
c_{12}^2 c_{21}^2 - 4 c_{03} c_{21}^3 - 4 c_{12}^3 c_{30} + 18 c_{03} c_{12} c_{21} c_{30} -  27 c_{03}^2 c_{30}^2.$$

\t
($\Pi^f$):
If $D\not=0$, then  we may assume $j^3f=x^2y\pm y^3$.
Then the $4$-jet can be transformed to
$$\Pi^f_{1}: \; xy^2\pm x^3 + \alpha x^3y + \beta y^4$$
by $\Psi$ with
$q_1=x$, $q_2=y$, $q_3=z$,
$p=1+\frac{1}{2}c_{22}x+\frac{1}{2}c_{13}y$.
Here  $\alpha=-c_{13}+c_{31}$ and $\beta=c_{04}$.
We may take either of $\alpha$ or $\beta$ to be $1$ unless both are $0$
(however it does not matter for latter discussion in Proposition \ref{BDE1}).

If $D=0$, then $j^3f\sim x^2y$ by some linear change.
Furthermore,  we may assume $c_{40},\, c_{04} (=\beta)\not=0$ generically.
Then by the same $\Psi$ above and a linear change, we have
$$\Pi^f_{2}: \;  xy^2+ x^4 \pm y^4+\alpha x^3y.$$

\section{Nets of asymptotic curves and their bifurcations}

\subsection{Binary differential equations}
In this section, we establish a relation between normal forms
in Theorem \ref{thm1} 
and local bifurcations of nets of asymptotic curves.
These curves are defined by a differential equation of particular form on the surface,
which had been studied in classical literatures (e.g. \cite{Tresse, Wil}). 
Here we take a modern approach from dynamical system. 
First we recall some basic notions in a general setting. 
A {\it binary differential equation} (BDE) in two variables $x, y$
has the form
$$F: \quad a(x,y)dy^2+2 b(x,y) dx dy + c(x,y) dx^2=0$$
with smooth functions $a, b, c$ of $x, y$. We can regard a BDE as
a map $\R^2 \to \R^3$ assigning $(x, y) \mapsto (a, b, c)$ and
consider the $C^\infty$-topology.
A BDE defines a pair of directions at each point $(x,y)$ in the plane
where $\delta(x,y):=b(x,y)^2-a(x,y)c(x,y)>0$.
Around a hyperbolic point, i.e. $\delta>0$,
the BDE locally defines two transverse foliations, thus
of our interest is the BDE around a point with $\delta=0$;
such points form the {\it discriminant} of the BDE.
At any point of the curve, the direction defined by BDE is unique,
and any solution of BDE passing the point generically has a cusp.
We say two germs of BDE's $F$ and $G$ are {\it topological equivalent} if
there is a local homeomorphism in the $xy$-plane
sending the integral curves of $F$ to those of $G$.

\begin{figure}
  \includegraphics[width=10cm]{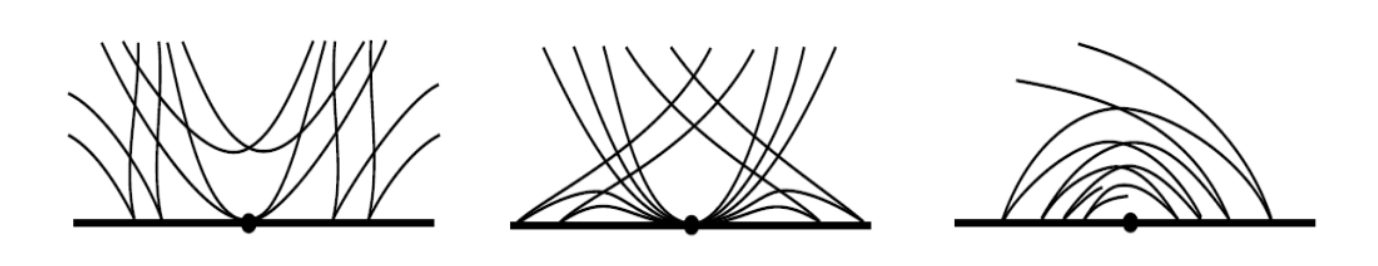}\\
  \caption{Folded saddle (left), node (center) and focus (right)}\label{fig1}
\end{figure}

One can separate BDE into two cases.
The first case occurs
when the functions $a, b, c$ do not vanish at the origin at once,
then the BDE is just an implicit differential equation.
The second case is that all the coefficients of BDE vanish at the origin.
Stable topological models of BDEs belong to the first case;
 it arises when the discriminant is smooth (or empty).
If the unique direction at any point of the discriminant is transverse to it
(i.e. integral curves form a family of cusps),
then the BDE is stable and smoothly equivalent to
$dy^2+xdx^2=0$, that was classically known in Cibrario \cite{CBR} and also Dara \cite{DR}.
If the unique direction is tangent to the discriminant,
then the corresponding point in the plane is called a \emph{folded singularity},
and  the BDE is stable and smoothly equivalent to
$dy^2+(-y+\lambda x^2)dx^2=0$
with $\lambda\neq 0,\frac{1}{16}$, as shown in Davydov \cite{Davydov, Davydov2}
-- there are three models of folded singularity,
see  Figure \ref{fig1}:
it is called a \emph{folded saddle} if $\lambda<0$, a \emph{folded node} if $0<\lambda<\frac{1}{16}$, and a \emph{folded focus} if $\frac{1}{16}<\lambda$.

In both cases,
the topological classification of
generic $1$ and $2$-parameter families of BDEs have been established
in Bruce-Fletcher-Tari \cite{BT1, BFT} and Tari \cite{Tari1, Tari2}, respectively.
We use those results below.
Besides, we need
a generic $3$-parameter family  of BDE studied in Oliver \cite{Oliver}.

\subsection{BDE of asymptotic curves}
Consider the surface locally given by Monge form $(x,y,f(x,y))$ around the origin.
Then asymptotic lines are defined by
$$f_{yy} dy^2 + 2 f_{xy} dx dy + f_{xx} dx^2=0,$$
which we call {\em asymptotic BDE} for short. 
The discriminant is the same as the parabolic curve. 
We first discuss bifurcations of asymptotic BDE around parabolic and flat umbilic points. 

Note that asymptotic BDE's form a particular class of general BDE's; 
the classification of asymptotic BDE is 
intimately related to the singularity type of the height function (i.e. the Monge form). 
The parabolic curve can be seen as the locus
where the height function has $A_{\geq2}$-singularities;
when the height function has a $A_3^\pm$-singularity,
the surface has an ordinary cusp of Gauss,
which corresponds to the type $\Pi^p_{c,1}$ --
in this case the asymptotic BDE has
a folded saddle singularity ($A_3^-$)
and has a folded node or focus singularity ($A_3^+$).
The transitions in $1$-parameter families occur generically
in three ways at the following singularities of the height function:
non-versal $A_3$,  $A_4$ and  $D_4$ (flat umbilic) \cite{BFT}.
For $2$-parameter families,
$A_3,$  $A_4$,  $A_5$ and $D_5$ singularities of the height function
generically appear.

Combining these results  with our Theorem \ref{thm1},
the following propositions are immediately obtained.

\begin{prop} \label{BDE0}
The following classes in Table \ref{main_table1}
correspond to structurally stable types of BDE given  in \cite{CBR, DR, Davydov}. 
\begin{itemize}
\item  $\Pi^p_{\I,k}\; (1 \le k\le 4):$ 
the parabolic curve is smooth
and the unique direction defined by $\delta=0$ is transverse to the curve; 
the BDE is smoothly equivalent to
$$dy^2+x\, dx^2=0.$$
\item  $\Pi^p_{c,1}, \; \Pi^p_{c,4},\; \Pi^p_{c,5}:$ 
the parabolic curve is smooth
and the unique direction defined by  $\delta=0$ is tangent to the curve, 
i.e. it is an ordinary/degenerate cusp of Gauss; 
the BDE is smoothly equivalent to
$$dy^2+(-y+\lambda x^2)dx^2=0$$
with $\lambda=6c_{40}-\frac{3}{2}\not=0$,
where $c_{40}$ is the coefficient of $x^4$ in the normal form.
\end{itemize}
\end{prop}

\proof The results follow from the comments in section 4.1 above.
If $j^4f= y^2+x^2y+c_{40}x^4$,
the $2$-jet of the asymptotic BDE is transformed to the above form (below)
via $x=\bar{x}$ and $y=-\frac{1}{2}\bar{x}^2-\bar{y}$.
\qed

\begin{rem}\label{remark_BDE0}{\rm
As $c_{40}=0$ for $\Pi^p_{c,4}$ and $\Pi^p_{c,5}$, we see $\lambda=-\frac{3}{2}<0$,
thus the BDE has a folded saddle at the origin. 
The {\it folded saddle-node bifurcation} (cf. Fig.2 in \cite{Tari1}) occurs at $\lambda=0$. 
That is the case of $c_{40}=\frac{1}{4}$,
that corresponds to the classes $\Pi^p_{c,k}\; (k=2,3)$ dealt below. 
However,  another exceptional value $\lambda=\frac{1}{16}$ does not relate to
the projection type, i.e.
the {\it  folded node-forcus bifurcation} of asymptotic BDE occurs
within the same class $\Pi^p_{c,1}$ (cf. Fig.3 in \cite{Tari1}). 
}
\end{rem}

\begin{prop}\label{BDE1}
The following classes  in Table \ref{main_table1}
correspond to some topological types of BDE with codimension $1$.
\begin{itemize}
\item $\Pi^p_{v,1}:$
the Monge form has
an $A_3$-singularity at the origin, at which
the parabolic curve has a Morse singularity;
the BDE is topologically equivalent to the non-versal $A_3^\pm$-transitions with
Morse type 1  in \cite{BFT}
$$dy^2+(\pm x^2 \pm y^2) dx^2=0.$$
\item $\Pi^p_{c,2}:$
the Monge form has
an $A_4$-singularity at the origin, at which
the parabolic curve is smooth ;
the BDE is topologically equivalent to the well-folded saddle-node type in \cite{BFT, Davydov2}
$$dy^2+(-y + x^3) dx^2=0,$$
provided the coefficient of $x^5$ in the normal form $c_{50}\neq 0$.
\item $\Pi^f_{1}:$
the Monge form has
a $D_4^\pm$-singularity at the origin, at which
the parabolic curve has a Morse singularity;
the BDE is topologically equivalent to the bifurcation of star/$1$-saddle types in \cite{BT1}
\begin{eqnarray*}
D_4^+: && ydy^2-2x dxdy -y dx^2=0 \;\; \mbox{\rm (star)};\\
D_4^-: &&ydy^2+2x dxdy +y dx^2=0 \;\; \mbox{\rm ($1$-saddle)}.
\end{eqnarray*}
\end{itemize}
\end{prop}

\proof
For $\Pi^p_{v,1}$,
our conditions  $c_{30}=c_{21}=0$ in  Table \ref{criteria_p}
appears in \cite[p.501]{BFT}, and also
$S\not=0$ if and only if the $2$-jet $j^2\delta(0)$ is non-degenrate,
hence the normal form follows from Theorem 2.7 (and Prop. 4.1)  in \cite{BFT}.
For $\Pi^p_{c,2}$,
$c_{30}=B=0$ and $A\not=0$ in Table \ref{criteria_p}
completely coincide with the condition for $A_4$-transition in  \cite[p.502]{BFT},
and the normal form of BDE is immediately obtained  (cf. \cite{Davydov2}).
For the class $\Pi^f_{1}$,
the asymptotic BDE is given in \cite[Cor. 5.3]{BT1}:
indeed, for the normal form of $\Pi^f_{1}$,
the parabolic curve is defined by
$3 x^2 -  y^2 + 18 \beta x y^2+ \cdots =0$,
hence it has a node at the origin for any parameters $\alpha, \beta$.
\qed

\begin{prop}\label{BDE2}
The following classes in Table \ref{main_table1}
correspond to some topological types of BDE with codimension $\ge 2$.
Below $c_{ij}$ are coefficients of the normal form of each class in Table \ref{main_table1}.

\begin{itemize}
\item  $\Pi^p_{v,2}:$
the Monge form has
an $A_3$-singularity at the origin,
at which the parabolic curve has a cusp singularity;
the BDE is topologically equivalent to the cusp type  in \cite{Tari1}
$$dy^2+(\pm x^2+y^3) dx^2=0,  $$
provided 
$\lambda_1:=\pm5c_{50}c_{31}^3+12c_{41}c_{31}^2\pm24c_{32}c_{31}+32c_{23}\neq  0
$.
\item  $\Pi^p_{v,3}:$
the Monge form has
an $A_4$-singularity at the origin,
at which the parabolic curve has a Morse singularity;
the BDE is topologically equivalent to the non-transvese Morse type  in \cite{Tari1}
$$dy^2+(xy+x^3) dx^2=0$$
provided $\lambda_2:=c_{31}\neq 0$.
\item  $\Pi^p_{c,3}:$
the Monge form has an $A_5$-singularity at the origin,
at which the parabolic curve is smooth;
the BDE is topologically equivalent to the folded degenerate elementary type  in \cite{Tari1}
$$dy^2+(-y \pm x^4) dx^2=0,$$
provided $\lambda_3:=c_{60}-\frac{1}{2}c_{41}\neq 0$.
\item $\Pi^f_{2}:$
the Monge form has
a $D_5$-singularity at the origin, at which
the parabolic curve has a cusp singularity;
the BDE is topologically equivalent to a cusp type 2 in \cite{Oliver}
$$x^2dx^2+2ydxdy+xdy^2=0.$$
\end{itemize}
Here $\lambda_i$ coincide moduli parameters 
appearing in the normal forms of jets of BDE studied in \cite{Tari1}. 
\end{prop}

\proof
For the first three classes,
it suffices to check the condition for BDE
in Proposition 4.1 of  \cite{Tari1}.
Indeed, for each of three classes,
our condition described in Table \ref{criteria_p} of \S \ref{criteria}
and additional conditions described above, $\lambda_1, \lambda_2, \lambda_3\not=0$,
are found in \cite[p.156]{Tari1}, and
the jets of BDE are given there
(see Theorem 1.1 of \cite{Tari1} for the normal forms).
For  $\Pi^f_{2}$,
as the $1$-jet of the asymptotic BDE is given by $j^1(a,b,c)(0)=(2x,2y,0)$ and
the parabolic curve is defined by $- 4 y^2+24 x^3 +\cdots =0$,
namely it has a cusp at the origin, the result follows from Theorem 3.4 in  \cite{Oliver}.
\qed

\begin{rem}{\rm
In Oliver \cite{Oliver},
BDE with the discriminant having a cusp are classified up to codimension $4$,
that generalizes a result for cusp types of codimension $2$ in \cite{Tari2}.
It is remarkable that the BDE for $\Pi^f_{2}$  in Proposition \ref{BDE2} is
one of such BDE with codimension $3$ in the space of all BDE,
while the stratum $\Pi^f_{2}$ has codimension $2$ in the space of all Monge forms 
as seen in Theorem \ref{thm1}. 
}
\end{rem}

\subsection{Bifurcations of parabolic and flecnodal curves}
Smooth and topological classifications of BDE do miss 
the geometry of flecnodal curves,  although it is quite fruitful.  
In part of his dissertation  \cite{UV2} Uribe-Vargas presented 14 types 
of generic $1$-parameter bifurcations of 
parabolic and flecnodal curves and of special elliptic points on the surface. 
That was done by a nice geometric argument using projective duality, 
however explicit normal forms are not given there. 
Thus we try to detect substrata in our classification of Monge forms 
which correspond to Uribe-Vargas' types. 

Each type defines an invariant open subset of one of our class of codimension $3$, 
or a proper invariant subset of a class of less codimension, which should be  
characterized by some particular relation of moduli parameters in the normal form. 
In fact we have already seen 8 types of degenerate parabolic and flat umbilical points in Uribe-Vargas' list 
as in all cases in Proposition \ref{BDE1} 
and an exceptional case in Proposition \ref{BDE0}: 
\begin{itemize}
\item 
$\Pi^p_{v,1}$ (4 types w.r.t. $\pm x^4$ and $\beta \lessgtr \frac{3}{8}\alpha^2$ in the normal form); 
\item 
$\Pi^p_{c,2}$ (the folded saddle-node type); 
\item 
$\Pi^f_1$ (star type and $1$-saddle type);  
\item
$\Pi^p_{c,4}$ (degenerate cusp of Gauss). 
\end{itemize}

\begin{rem}
{\rm 
The corresponding strata in jet space have codimension $3$. 
Take a generic family of surfaces passing through each stratum. 
Around a point of class $\Pi^p_{v,1}$, 
Morse bifurcations of the parabolic curve happen, 
while bifurcations of the flecnodal curve  
can not be analyzed from the classification of BDE in \cite{BFT}. 
Such bifurcations have been determined completely in \cite{UV2} 
and can also be checked by direct computations from our Monge forms. 
For $\Pi^p_{c,2}$, 
two folded saddle points are cancelled at this point, 
called in  \cite[Prop. 4.2]{BFT} the folded saddle-node bifurcation. 
For $\Pi^f_1$, there are two types of degenerate flat umbilic points. 
For $\Pi^p_{c,4}$, 
parabolic and flecnodal curves meet tangentially, and 
a butterfly point comes across at the cusp of Gauss point. 
As an additional remark, for $\Pi^p_{c,1}$ with $\lambda=\frac{1}{16}$, 
the folded node-forcus bifurcation of BDE appears as mentioned in Remark \ref{remark_BDE0}, 
but the parabolic and flecnodal curves do not bifurcate, 
so this case is not included. 

}
\end{rem}

The rest are divided into 5 types occurred at hyperbolic points and 1 type at an elliptic point; 
the normal forms for degenerate hyperbolic points are given below: 
\begin{itemize}
\item 
$\Pi^h_{3,6}$ (degenerate butterfly type); 
\item 
$\Pi^h_{4,5}$ (swallowtail+butterfly);  
\item 
$\Pi^h_{3,5}$ with $\alpha=0$  (Morse bifurcation; `lips' and `bec-\`a-bec'); 
\item
$\Pi^h_{4,4}$ with $\alpha\beta=16$ (tacnode bifurcation). 
\end{itemize}

In  Table \ref{main_table2} (hyperbolic Monge forms), there are two strata of codimension $3$. 
For the normal form of $\Pi^h_{3,6}$, the flecnodal curve is smooth, and 
the $y$-axis is tangent to the surface with $6$-point contact. 
Through a generic bifurcation, two butterfly points ($\Pi^h_{3,5}$) are cancelled (or created) at this point. 
For $\Pi^h_{4,5}$, 
there are two irreducible components of the flecnodal curve, 
one of which has a butterfly at the origin. 
It bifurcates into one point of $\Pi^h_{4,4}$ and one butterfly point. 

In strata of codimension $2$, 
generic bifurcations of flecnodal curves may occur 
at some particular values of parameters appearing in the normal forms. 
First consider the class $\Pi^h_{3,k}\; (k\ge 4)$: 
$f(x,y)=xy+x^3+\sum_{i+j\ge 4} c_{ij}\, x^i y^j$ ($\alpha=c_{13}$ in Table \ref{main_table2}). 
We are now viewing the surface along lines close to the $y$-axis, 
so take parallel projection $\varphi: (x,y) \mapsto (x- uy, f(x,y)-v y)$. 
The flecnodal curve is just the locus of singular points $(x,y)$ at which 
the projection for some $(u,v)$ has the swallowtail singularities or more degenerate ones. 
Put $\lambda:=\det d\varphi$ and 
$\eta:=u \frac{\rd}{\rd x}+ \frac{\rd}{\rd y}$ 
so that $\eta$ spans the kernel of $d\varphi$ over $\lambda=0$. 
Then the locus is defined to be the image via projection $(x,y,u,v) \mapsto (x,y)$ 
of curves given by three equations $\lambda=\eta \lambda=\eta\eta \lambda=0$ (cf. \cite{Kabata}). 
By $\lambda=0$ we eliminate $v$, and by $\eta\lambda=0$, 
we solve $u$ by implicit function theorem: 
$u=-c_{22}x^2-3c_{13}xy-6c_{04}y^2+o(2)$. 
Substitute them into $\eta\eta\lambda=0$; 
we have the expansion of the equation of flecnodal curve: 
$$c_{13} x + 4 c_{04} y +  c_{23} x^2 + 4 c_{14} x y + 10 c_{05} y^2+h.o.t.=0.$$ 
The curve is not smooth at the origin when $c_{13}=c_{04}=0$, and 
it has a Morse singularity  if $5c_{23}c_{05}-2c_{14}^2\not=0$. 
Thus, we conclude that 
Morse bifurcations of flecnodal curves occur 
at a point of the class $\Pi^h_{3,5}$ with $\alpha=0$ and generic homogeneous terms of order $5$ 
(those are regarded as normal forms for  `lips' and `bec-\`a-bec' in \cite{UV2}). 

Next, consider the class $\Pi^h_{4,4}$: 
$f(x,y)=xy+x^4+y^4+\alpha xy^3+\beta x^3y+o(4)$. 
In entirely the same way as seen above, 
for each of projections 
$\varphi_L: (x,y) \mapsto (x- uy, f(x,y)-v y)$ and 
$\varphi_R: (x,y) \mapsto (y- ux, f(x,y)-v x)$, 
we take $\lambda$ and $\eta$, and solve $\lambda=\eta \lambda=\eta\eta \lambda=0$. 
We then obtain the equation of each component of the flecnodal curve: 
$$\alpha x + 4 y+ h.o.t=0, \quad 4x+\beta y+ h.o.t=0.$$ 
Thus 
the tacnode bifurcation (tangency of two components) occurs 
at a point of the class  $\Pi^h_{4,4}$ with $\alpha\beta=16$.  
The condition exactly coincides with the one in Table 3.2 in Landis \cite{Landis}. 
It is also observed that 
in the elliptic domain, 
a similar bifurcation occurs at the class $\Pi^e_{4,4}$ with a particular value.


\end{document}